\documentclass{amsart}

\usepackage{amsmath,amsthm,amssymb,txfonts,bm,xcolor}
\usepackage[utf8]{inputenc}

\theoremstyle{definition}
\newtheorem{definition}{Definition}
\newtheorem{theorem}[definition]{Theorem}
\newtheorem{lemma}[definition]{Lemma}

\newtheorem{remark}[definition]{Remark}

\def\p#1#2#3#4{{}_{2}\phi_{1}\biggl(\genfrac..{0pt}{}{{#1},\,{#2}}{#3};q, #4\biggr)}
\def\opo#1#2#3{{}_{1}\phi_{1}\biggl(\genfrac..{0pt}{}{#1}{#2};q, #3\biggr)}
\def\zpo#1#2{{}_{0}\phi_{1}\biggl(\genfrac..{0pt}{}{\text{--}}{#1};q, #2\biggr)}

\def\tpt#1#2#3#4#5#6{{}_{3}\phi_{2}\biggl(\genfrac..{0pt}{}{{#1},\,{#2},\,{#3}}{{#4},\,{#5}};q, #6\biggr)}
\def\c#1#2#3#4#5#6#7#8{\biggl(\genfrac..{0pt}{}{{#1},\,{#2}}{#3};{#4}\,;\genfrac..{0pt}{}{{#5},\,{#6}}{#7};q, #8\biggr)}

\makeatletter
\@addtoreset{equation}{section}

\makeatother

\title[Three-Term Recurrence Relations for Confluent Basic Hypergeometric Series]{Three-Term Recurrence Relations for Confluent Basic Hypergeometric Series with Applications to $q$-Bessel Functions}
\author{Yuka Yamaguchi}
\date{\today}

\begin{document}

\begin{abstract}
We establish three-term recurrence relations for the ${}_1\phi_1$ and ${}_0\phi_1$ basic hypergeometric series involving multiplicative shifts of the parameters and the variable by integer powers of $q$. The coefficients of these recurrence relations are shown to be uniquely determined by the shift indices and are given explicitly in terms of rational functions. These recurrence relations arise as confluent limits of previously established recurrence relations for the ${}_2\phi_1$ basic hypergeometric series.

As an application, we derive three-term recurrence relations for Jackson's second and third $q$-Bessel functions. These recurrence relations involve additive shifts in the order and multiplicative $q$-shifts in the variable, and their coefficients include the known $q$-Lommel polynomials as special cases.


\end{abstract}

\maketitle

\section{Introduction}
For the ${}_2\phi_1$ basic hypergeometric series, a general three-term recurrence relation involving multiplicative shifts of the parameters and the variable by integer powers of $q$ was established in \cite{Y1}. In that work, the coefficients are shown to be uniquely determined by the shift indices and are given explicitly as rational functions. These results provide a natural starting point for investigating corresponding recurrence relations for the ${}_1\phi_1$ and ${}_0\phi_1$ basic hypergeometric series through confluent limit procedures.

In this paper, we establish three-term recurrence relations for the ${}_1\phi_1$ basic hypergeometric series involving multiplicative shifts by integer powers of $q$, prove that the coefficients are uniquely determined by the shift indices, and obtain explicit expressions for them in terms of rational functions. We then derive the corresponding recurrence relation for the ${}_0\phi_1$ basic hypergeometric series by taking a confluent limit. Finally, the recurrence relations for the ${}_1\phi_1$ and ${}_0\phi_1$ basic hypergeometric series are applied to Jackson's third and second $q$-Bessel functions, respectively. These recurrence relations involve additive shifts in the order and multiplicative $q$-shifts in the variable, and their coefficients include the known $q$-Lommel polynomials \cite{Ismail, KS} as special cases. 

Let $(a; q)_{j}$ and $(a; q)_{\infty}$ denote the $q$-shifted factorials defined by 
\begin{align*}
(a; q)_{j} &:= 
\begin{cases}
1, & j=0, \\
(1-a)(1-aq) \cdots (1-aq^{j-1}), & j \geq 1, 
\end{cases}\\
(a; q)_{-j} &:= \frac{1}{(aq^{-j};q)_{j}} = \frac{1}{(1-aq^{-1})(1-aq^{-2}) \cdots (1-aq^{-j})}, \quad j \geq 1, \\
(a; q)_{\infty} &:= \prod_{j = 0}^{\infty} (1 - a q^{j}). 
\end{align*}
The ${}_r\phi_s$ basic hypergeometric series is defined by 
\begin{align*}
{}_{r}\phi_{s}\biggl(\genfrac..{0pt}{}{a_{1},\,a_{2}, \dotsc, a_{r}}{b_{1}, \dotsc, b_{s}}; q, x\biggr) 
&= {}_{r}\phi_{s} (a_{1}, a_{2}, \ldots, a_{r}; b_{1}, \ldots, b_{s}; q, x) \\
&:= \sum_{j = 0}^{\infty} 
\frac{(a_{1}; q)_{j} (a_{2}; q)_{j} \dotsm (a_{r}; q)_{j}}{(q; q)_{j} (b_{1}; q)_{j} \dotsm (b_{s}; q)_{j}} 
\left\{(-1)^j q^{\binom{j}{2}} \right\}^{1+s-r} x^{j}, 
\end{align*}
where $\binom{j}{2} := j(j-1)/2$. This definition follows \cite[(1.2.22)]{GR}. 
It is assumed that none of the lower parameters $b_{1}, \dotsc, b_{s}$ equals $1$ or any negative integer power of $q$. 
In this paper, we assume $|q|<1$.

For any integers $k,l,m,n$, the following three-term recurrence relation for ${}_2\phi_1$ is considered in \cite[(1.3)]{Y1}: 
\begin{align}\label{3trr_2phi1}
\p{a q^{k}}{b q^{l}}{c q^{m}}{x q^{n}} 
= Q \cdot \p{a q}{b q}{c q}{x} + R \cdot \p{a}{b}{c}{x}. 
\end{align}
Here $Q$ and $R$ are rational functions in $a, b, c, x,$ and $q$. 
In \cite[Theorem~1]{Y1}, it is proved that the pair $(Q, R)$ satisfying $(\ref{3trr_2phi1})$ is uniquely determined by $(k, l, m, n)$, and explicit expressions for $Q$ and $R$ are given. 
For the special cases with $n=0$, see \cite[Lemma~1, Theorem~2]{Y0}. 

Although the assumption 
\begin{align*}
a, b, c, \frac{a}{b}, \frac{c}{a}, \frac{c}{b} \notin q^{\mathbb{Z}} \cup\{0 \}
\end{align*}
is imposed in \cite[(1.4)]{Y1}, whose results are used in this paper, 
we do not explicitly restrict the range of the parameters here. 
Throughout this paper, the parameters are understood to take values in the region where the series converge and the rational functions are defined; this restriction is justified by analytic continuation.

In this paper, we first consider two formulations of the three-term recurrence relation for ${}_1\phi_1$. 
We begin with the following formulation. 
For any integers $k,m,n$, there exist rational functions
\[
S = S(k,m,n; a,c,x; q), \qquad
T = T(k,m,n; a,c,x; q)
\]
in $a,c,x,$ and $q$ satisfying 
\begin{align}\label{3trr_1phi1}
\opo{a q^k}{c q^m}{x q^n}
=
S \cdot\opo{a q}{c q}{x q}
+
T \cdot\opo{a}{c}{x}.
\end{align}

\begin{theorem}\label{thm:1}
The pair $(S, T)$ satisfying $(\ref{3trr_1phi1})$ is uniquely determined by $(k,m,n)$. 
Consequently, $S$ and $T$ satisfy the relation
\begin{align}\label{relation_S,T}
T(k,m,n; a,c,x; q) 
= - \frac{(1-c)(1-cq)}{(1-aq)(c-axq)xq}\, 
S(k-1,m-1,n-1; aq,cq,xq; q). 
\end{align}
Moreover, $S$ is given explicitly by 
\begin{align}\label{expression_S}
S(k,m,n; a,c,x; q) 
= \frac{(1-a)(c-ax)}{(q-c)(1-c)} 
\frac{x^{1-\max\{m, 0\}}}{(ax/c;q)_{\max\{k-m+n, 0\}}}
P^{(1,1)}(k,m,n; a,c,x; q), 
\end{align}
where $P^{(1,1)}(k,m,n; a,c,x; q)$ is a polynomial in $x$ of degree at most 
\[
\max \{|k|, |m-k|, |n|, |m-n|\} - 1,
\]
whose explicit expression is given in $\eqref{P^(1,1)}$. 
\end{theorem}

We next consider an alternative formulation of the recurrence relation, 
again for arbitrary integers $k,m,n$, in which
the shift in the variable of ${}_1\phi_1$ is arranged differently. 

\begin{align}\label{3trr_1phi1_2}
\opo{a q^k}{c q^m}{x q^n}
=
\tilde{S} \cdot\opo{a q}{c q}{x}
+
\tilde{T} \cdot\opo{a}{c}{x}, 
\end{align}
where 
\[
\tilde{S} = \tilde{S}(k,m,n; a,c,x; q), \qquad
\tilde{T} = \tilde{T}(k,m,n; a,c,x; q)
\]
are rational functions in $a,c,x,$ and $q$. 

\begin{theorem}\label{thm:2}
The pair $(\tilde{S}, \tilde{T})$ satisfying $(\ref{3trr_1phi1_2})$ is uniquely determined by $(k,m,n)$. 
Consequently, they satisfy the relation
\begin{align}\label{relation_tS,tT}
\tilde{T}(k,m,n; a,c,x; q) 
= - \frac{(1-c)(1-cq)q}{(1-aq)x}\,  
\tilde{S}(k-1,m-1,n; aq,cq,x; q). 
\end{align}
Moreover, $\tilde{S}$ is given explicitly by 
\begin{align}\label{expression_tS}
\tilde{S}(k,m,n; a,c,x; q) 
= \frac{1-a}{(q-c)(1-c)} 
\frac{x^{1-\max\{m, 0\}}}{(ax/c;q)_{\max\{k-m+n, 0\}}}
P^{(1,1)}(k,m,n; a,c,x; q), 
\end{align}
where $P^{(1,1)}(k,m,n; a,c,x; q)$ is a polynomial in $x$ of degree at most 
\[
\max \{|k|, |m-k|, |n|, |m-n|\} - 1,
\]
whose explicit expression is given in $\eqref{P^(1,1)}$. 
\end{theorem}

The recurrence relation \eqref{3trr_1phi1} arises naturally as a
confluent limit of \eqref{3trr_2phi1}. 
The formulation \eqref{3trr_1phi1_2} is particularly suited to applications to
$q$-Bessel functions. 

We now turn to the three-term recurrence relation for ${}_0\phi_1$
obtained as a confluent limit of \eqref{3trr_1phi1_2}. 
For any integers $m,n$, there exist rational functions
\[
U = U(m,n; c,x; q), \qquad
V = V(m,n; c,x; q)
\]
in $c,x,$ and $q$ satisfying
\begin{align}\label{3trr_0phi1}
\zpo{c q^m}{x q^n}
=
U \cdot\zpo{c q}{x q}
+
V \cdot\zpo{c}{x}.
\end{align}
Let $\lfloor a \rfloor$ denote the greatest integer not exceeding $a$. 

\begin{theorem}\label{thm:3}
The pair $(U, V)$ satisfying $(\ref{3trr_0phi1})$ is uniquely determined by $(m,n)$. 
Consequently, they satisfy the relation
\begin{align}\label{relation_U,V}
V(m,n; c,x; q) 
= \frac{(1-c)(1-cq)}{x} 
U(m-1,n-1; cq,xq; q). 
\end{align}
Moreover, $U$ is given explicitly by 
\begin{align}\label{expression_U}
U(m,n; c,x; q) 
= -\frac{1}{(q-c)(1-c)} 
\frac{x^{1-\max\{m, 0\}}}{(x/c;q)_{\max\{n-m, 0\}}}
P^{(0,1)}(m,n; c,x; q), 
\end{align}
where $P^{(0,1)}(m,n; c,x; q)$ is a polynomial in $x$ of degree at most 
\[
\left\lfloor \frac{n+1}{2} \right\rfloor - \min\{m, n, n-m, 0\} - 1, 
\]
whose explicit expression is given in $\eqref{P^(0,1)}$.  
\end{theorem}

We apply the recurrence relations obtained above to Jackson's $q$-Bessel functions. 

Jackson introduced the following three $q$-analogues of Bessel functions
(see, for example, \cite[p.~30]{GR} and \cite[p.~354]{Ismail_book}):
\begin{align*}
J^{(1)}_\nu (x; q) 
&:= \frac{(q^{\nu + 1};q)_\infty}{(q;q)_\infty} \left(\frac{x}{2} \right)^\nu 
\p{0}{0}{q^{\nu+1}}{-\frac{x^2}{4}}, 
\qquad |x|<2, \\
J^{(2)}_\nu (x; q) 
&:= \frac{(q^{\nu + 1};q)_\infty}{(q;q)_\infty} \left(\frac{x}{2} \right)^\nu 
\zpo{q^{\nu+1}}{-\frac{x^2 q^{\nu+1}}{4}}, 
\qquad x \in \mathbb{C}, \\
J^{(3)}_\nu (x; q) 
&:= \frac{(q^{\nu + 1};q)_\infty}{(q;q)_\infty} \left(\frac{x}{2} \right)^\nu 
\opo{0}{q^{\nu+1}}{\frac{x^2q}{4}}, 
\qquad x \in \mathbb{C}. 
\end{align*}
Here $0 < q < 1$. 
We have 
\begin{align*}
\lim_{q \to 1} J^{(k)}_\nu ((1-q)x; q) = J_\nu (x), 
\qquad k = 1,2,3, 
\end{align*}
where $J_\nu (x)$ denotes the Bessel function of the first kind of order $\nu$. 

Using the above representations, we derive the three-term recurrence
relations for Jackson's third and second $q$-Bessel functions directly from the
recurrence relations for ${}_1\phi_1$ and ${}_0\phi_1$ established above.
Since
\[
J^{(1)}_\nu (x; q)
= \frac{J^{(2)}_\nu (x; q)}{(-x^2/4; q)_\infty},
\qquad |x|<2,
\]
(see \cite[(14.1.18)]{Ismail_book}),
the three-term recurrence relation for $J^{(1)}_\nu (x; q)$
follows directly from that of $J^{(2)}_\nu (x; q)$ and is therefore omitted.

We first apply the recurrence relation \eqref{3trr_1phi1_2} 
to Jackson's third $q$-Bessel function. 

\begin{theorem}\label{thm:4}
For any integers $m,n$, we have
\begin{align}\label{3trr_J^(3)}
q^{-\frac{(\nu+m)n}{2}} J^{(3)}_{\nu+m}(xq^{\frac{n}{2}};q) 
= R^{(3)}_{m,n,\nu}(x;q)\, J^{(3)}_\nu (x;q) 
- R^{(3)}_{m-1,n,\nu+1}(x;q)\, J^{(3)}_{\nu-1}(x;q), 
\end{align}
where 
\begin{align*}
R^{(3)}_{m,n,\nu} (x;q)
:= \frac{(q^{\nu+m+1};q)_\infty}{(q^{\nu-1};q)_\infty} 
q^{-\max\{m+1,0 \}} \left(\frac{x}{2} \right)^{\min\{m+2, -m \}} 
P^{(3)}_{m+1,n,\nu}(x;q), 
\end{align*}
and $P^{(3)}_{m,n,\nu}(x;q)$ is a polynomial in $x$ of degree 
at most $\max\{|n|, |m-n| \} - 1$, 
whose explicit expression is given in \eqref{P^(3)}. 
\end{theorem}

When $m \ge 1$ and $n=0$, after the change of the variable $x \mapsto 2x$, \eqref{3trr_J^(3)} coincides with the three-term recurrence relation considered in \cite[(4.12)]{KS}, 
and the Laurent polynomial $R^{(3)}_{m,n,\nu}(2x;q)$ in $x$ coincides with the $q$-Lommel polynomial $R_{m,\nu}(x;q)$ in \cite[(4.23)]{KS}. 
See Remark~\ref{remark_KS} for details.

We next apply the recurrence relation \eqref{3trr_0phi1}
to Jackson's second $q$-Bessel function. 

\begin{theorem}\label{thm:5}
For any integers $m,n$, we have
\begin{align}\label{3trr_J^(2)}
q^{-\frac{(\nu+m)n}{2} + m \nu + \binom{m}{2}} J^{(2)}_{\nu+m}(xq^{\frac{n}{2}};q) 
= R^{(2)}_{m,n,\nu}(x;q)\, J^{(2)}_\nu (x;q) 
- R^{(2)}_{m-1,n,\nu+1}(x;q)\, J^{(2)}_{\nu-1}(x;q), 
\end{align}
where 
\begin{align*}
R^{(2)}_{m,n,\nu} (x;q)
&:= \frac{(q^{\nu+m+1};q)_\infty}{(q^{\nu-1};q)_\infty (-x^2/4;q)_{\max\{n, 0\}}} \\
&\quad \times (-1)^{\max \{m+1, 0\}} q^{\min\{m+1,0 \} \nu + \binom{m}{2} - 1} \left(\frac{x}{2} \right)^{\min\{m+2, -m \}} 
P^{(2)}_{m+1,n,\nu}(x;q), 
\end{align*}
and $P^{(2)}_{m,n,\nu}(x;q)$ is a polynomial in $x$ of degree 
at most $\left\lfloor \frac{m+n+1}{2} \right\rfloor - \min\{m, n, m+n, 0 \} - 1$, 
whose explicit expression is given in \eqref{P^(2)}. 
\end{theorem}

When $m \ge 1$ and $n=0$, \eqref{3trr_J^(2)} coincides with the three-term recurrence relation considered in \cite[(1.19)]{Ismail}, 
and the rational function $R^{(2)}_{m,n,\nu}(x;q)$ in $x$ coincides with the $q$-Lommel polynomial $R_{m,\nu}(x;q)$ in \cite[(2.10)]{Ismail}. 
See Remark~\ref{remark_Ismail} for details.

\section{Explicit expressions for the polynomials}
In this section, we present explicit expressions for the polynomials appearing in Theorems~\ref{thm:1}--\ref{thm:5}. 
We note that 
\[
\frac{1}{(q;q)_j} = (q^{j+1};q)_{-j} = 0, \quad j<0. 
\]

\begin{definition}
Let $P^{(1,1)}(k,m,n; a,c,x; q)$ be the polynomial in $x$ defined by
\begin{align}\label{P^(1,1)}
P^{(1,1)}(k,m,n; a,c,x; q) 
:= 
\begin{cases}
\displaystyle\sum_{j=0}^{d} \left(A^{(1,1)}_{j+\min\{m, 0\}} - B^{(1,1)}_{j-\max\{m, 0\}} \right)x^j, 
& k-m+n \geq 0, \\
\displaystyle\sum_{j=0}^{d} \left(\tilde{A}^{(1,1)}_{j+\min\{m, 0\}} - \tilde{B}^{(1,1)}_{j-\max\{m, 0\}} \right)x^j, 
& k-m+n < 0, 
\end{cases}
\end{align}
where $d := \max\{k-m+n, 0\} + \max\{m,0\} - \min\{k,n\} - 1 = \max \{|k|, |m-k|, |n|, |m-n| \} - 1$, 
\begin{align*}
A^{(1,1)}_j 
&:= \frac{(aq/c;q)_{k-m} (c;q)_m (c/q;q)_{m-j}}{(q;q)_j (a;q)_{k-j}} 
c^{-n} q^{(j-m+1)(n-1)+1} 
\tpt{q^{-j}}{cq^{m-j-1}}{a}{c}{aq^{k-j}}{q^{j-n+1}}, \\
B^{(1,1)}_j 
&:= \frac{(aq/c;q)_{j}}{(q;q)_j (q^2/c;q)_{j}} 
(-1)^j c^{-j} q^{\binom{j+1}{2}} 
\tpt{q^{-j}}{cq^{-j-1}}{cq^{m-k}/a}{cq^m}{cq^{-j}/a}{q^{k-m+n+1}}, \\
\tilde{A}^{(1,1)}_j 
&:= \frac{(aq/c;q)_{j+k-m} (c;q)_m (c/q;q)_{m-j}}{(q;q)_j (a;q)_{k}} 
c^{-n} q^{(j-m+1)(n-1)+1} \\
&\qquad \times \tpt{q^{-j}}{cq^{m-j-1}}{c/a}{c}{cq^{m-k-j}/a}{q^{m-k-n+1}}, \\
\tilde{B}^{(1,1)}_j 
&:= \frac{(q/a;q)_{j}}{(q;q)_j (q^2/c;q)_{j}} 
a^j c^{-j} 
\tpt{q^{-j}}{cq^{-j-1}}{aq^k}{cq^m}{aq^{-j}}{q^{j+n+1}}. 
\end{align*}
\end{definition}
\begin{remark}\label{remark:P^(1,1)}
When $j<0$, we have 
$A^{(1,1)}_j = B^{(1,1)}_j = \tilde{A}^{(1,1)}_j = \tilde{B}^{(1,1)}_j = 0$. 
\end{remark}

\begin{definition}
Let $P^{(0,1)}(m,n; c,x; q)$ be the polynomial in $x$ defined by 
\begin{align}\label{P^(0,1)}
P^{(0,1)}(m,n; c,x; q) 
:= 
\begin{cases}
\displaystyle\sum_{j=0}^{e} \left(A^{(0,1)}_{j+\min\{m, 0\}} - B^{(0,1)}_{j-\max\{m, 0\}} \right)x^j, 
& m \leq n, \\
\displaystyle\sum_{j=0}^{e} \left(\tilde{A}^{(0,1)}_{j+\min\{m, 0\}} - \tilde{B}^{(0,1)}_{j-\max\{m, 0\}} \right)x^j, 
& m > n, 
\end{cases}
\end{align}
where $e := \left\lfloor \frac{n+1}{2} \right\rfloor - \min\{m, n, n-m, 0\} - 1$, 
\begin{align*}
A^{(0,1)}_j 
&:= \frac{(c;q)_m}{(q;q)_j (q^2/c;q)_{j-m}} 
c^{2m-n-j} q^{(j-m+1)(n-m)} 
\p{q^{-j}}{cq^{m-j-1}}{c}{q^{2j-n+1}}, \\
B^{(0,1)}_j 
&:= \frac{1}{(q;q)_j (q^2/c;q)_{j}} 
c^{-2j} q^{j(j+1)} 
\p{q^{-j}}{cq^{-j-1}}{cq^m}{q^{n-m+1}}, \\
\tilde{A}^{(0,1)}_j 
&:= \frac{(c;q)_m}{(q;q)_j (q^2/c;q)_{j-m}} 
c^{2m-n-2j} q^{(j-m+1)(j-m+n)} 
\p{q^{-j}}{cq^{m-j-1}}{c}{q^{m-n+1}}, \\
\tilde{B}^{(0,1)}_j 
&:= \frac{1}{(q;q)_j (q^2/c;q)_{j}} c^{-j} 
\p{q^{-j}}{cq^{-j-1}}{cq^m}{q^{2j+n+1}}. 
\end{align*}
\end{definition}
\begin{remark}\label{remark:P^(０,1)}
When $j<0$, we have 
$A^{(0,1)}_j = B^{(0,1)}_j = \tilde{A}^{(0,1)}_j = \tilde{B}^{(0,1)}_j = 0$. 
\end{remark}

\begin{definition}
Let $P^{(3)}_{m,n,\nu}(x;q)$ be the polynomial in $x$ defined by 
\begin{align}\label{P^(3)}
P^{(3)}_{m,n,\nu}(x;q) 
:= \sum_{j=0}^{f} \left(\frac{x}{2} \right)^{2j} q^j 
\left(A^{(3)}_{m,n,\nu,j+\min\{m,0\}} - B^{(3)}_{m,n,\nu,j-\max\{m,0\}} \right), 
\end{align}
where $f := \max\{|n|, |m-n| \} - 1$, 
\begin{align*}
A^{(3)}_{m,n,\nu,j} 
&:= \frac{(q^\nu;q)_m (q^{\nu-1};q)_{m-j}}{(q;q)_j} 
q^{(j-m+1)(n-1) + 1 - n \nu} 
\p{q^{-j}}{q^{\nu+m-j-1}}{q^\nu}{q^{j-n+1}}, \\
B^{(3)}_{m,n,\nu,j} 
&:= \frac{1}{(q;q)_j (q^{2-\nu};q)_j} 
(-1)^j q^{\binom{j+1}{2} - j \nu} 
\p{q^{-j}}{q^{\nu-j-1}}{q^{\nu+m}}{q^{j+n+1}}. 
\end{align*}
\end{definition}
\begin{remark}\label{remark:P^(3)}
When $j<0$, we have 
$A^{(3)}_{m,n,\nu,j} = B^{(3)}_{m,n,\nu,j} = 0$.  
\end{remark}
\begin{remark}\label{remark_KS}
When $m \ge 1$ and $n=0$, by Remark~\ref{remark:P^(3)}, 
and after a straightforward simplification, we have 
\begin{align*}
R^{(3)}_{m,0,\nu} (2x;q) 
= \sum_{j=0}^{m} x^{2j-m} 
\frac{(q^{j+1};q)_\infty (q^\nu;q)_\infty}{(q;q)_\infty (q^{\nu+m-j};q)_\infty} 
\p{q^{-j}}{q^{\nu+m-j}}{q^\nu}{q^{j+1}}.   
\end{align*}
This coincides with $R_{m,\nu}(x;q)$ in \cite[(4.23)]{KS}.
\end{remark}

\begin{definition}
Let $P^{(2)}_{m,n,\nu}(x;q)$ be the polynomial in $x$ defined by 
\begin{align}\label{P^(2)}
P^{(2)}_{m,n,\nu}(x;q) 
:= 
\begin{cases}
\displaystyle 
\sum_{j=0}^{g} \left(\frac{x}{2} \right)^{2j} (-1)^j q^{j \nu}
\left(A^{(2)}_{m,n,\nu,j+\min\{m,0\}} - B^{(2)}_{m,n,\nu,j-\max\{m,0\}} \right), & n \geq 0, \\
\displaystyle
\sum_{j=0}^{g} \left(\frac{x}{2} \right)^{2j} (-1)^j q^{j \nu} 
\left(\tilde{A}^{(2)}_{m,n,\nu,j+\min\{m,0\}} - \tilde{B}^{(2)}_{m,n,\nu,j-\max\{m,0\}} \right), & n < 0,
\end{cases}
\end{align}
where $g := \left\lfloor \frac{m+n+1}{2} \right\rfloor - \min\{m, n, m+n, 0\} - 1$, 
\begin{align*}
A^{(2)}_{m,n,\nu,j} 
&:= \frac{(q^\nu;q)_m}{(q;q)_j (q^{2-\nu};q)_{j-m}} 
q^{(j-m+1)n - (j-m+n) \nu} 
\p{q^{-j}}{q^{\nu+m-j-1}}{q^\nu}{q^{2j-m-n+1}}, \\
B^{(2)}_{m,n,\nu,j} 
&:= \frac{1}{(q;q)_j (q^{2-\nu};q)_j} 
q^{j(j+1) - 2j \nu} 
\p{q^{-j}}{q^{\nu-j-1}}{q^{\nu+m}}{q^{n+1}}, \\
\tilde{A}^{(2)}_{m,n,\nu,j} 
&:= \frac{(q^\nu;q)_m}{(q;q)_j (q^{2-\nu};q)_{j-m}} 
q^{(j-m+1)(j+n) - (2j-m+n) \nu} 
\p{q^{-j}}{q^{\nu+m-j-1}}{q^\nu}{q^{1-n}}, \\
\tilde{B}^{(2)}_{m,n,\nu,j} 
&:= \frac{1}{(q;q)_j (q^{2-\nu};q)_j} 
q^{- j \nu} 
\p{q^{-j}}{q^{\nu-j-1}}{q^{\nu+m}}{q^{2j+m+n+1}}. 
\end{align*}
\end{definition}
\begin{remark}\label{remark:P^(2)}
When $j<0$, we have 
$A^{(2)}_{m,n,\nu,j} = B^{(2)}_{m,n,\nu,j} 
= \tilde{A}^{(2)}_{m,n,\nu,j} = \tilde{B}^{(2)}_{m,n,\nu,j} = 0$. 
\end{remark}
\begin{remark}\label{remark_Ismail}
When $m \ge 1$ and $n=0$, by Remark~\ref{remark:P^(2)}, 
and after a straightforward simplification 
using the $q$-Chu-Vandermonde summation (see, e.g., \cite[(1.5.2)]{GR})
\begin{align*}
\p{q^{-j}}{b}{c}{\frac{cq^j}{b}} = \frac{(c/b;q)_j}{(c;q)_j}, 
\end{align*}
we obtain 
\begin{align*}
R^{(2)}_{m,0,\nu} (x;q) 
= \sum_{j=0}^{\lfloor m/2 \rfloor} \left(\frac{x}{2}\right)^{2j-m} 
\frac{(-1)^j (q^\nu;q)_{m-j} (q;q)_{m-j}}{(q;q)_j (q^{\nu};q)_j (q;q)_{m-2j}} 
q^{j(\nu+j-1)}. 
\end{align*}
This coincides with $R_{m,\nu}(x;q)$ in \cite[(2.10)]{Ismail}.
\end{remark}

\section{Three-term recurrence relation for ${}_1\phi_1$}\label{proof_thm1}
In this section, we prove Theorem~\ref{thm:1}. 

First, we prove the uniqueness of the pair $(S, T)$ satisfying \eqref{3trr_1phi1}. 
To prove this, we use the following three-term recurrence relation for ${}_2\phi_1$ obtained from \cite[Theorem~1]{Y1}: 
\begin{align*}
\p{aq}{b}{cq}{x} 
= \frac{(1-b)(c-abx)}{c-b}\, \p{aq}{bq}{cq}{x} 
- \frac{b(1-c)}{c-b}\, \p{a}{b}{c}{x}. 
\end{align*}
Replacing $x$ with $x/b$ and letting $b \to \infty$, we obtain 
\begin{align}\label{3trr_1phi1_(1,1,0)}
\opo{aq}{cq}{x} 
= (c-ax)\, \opo{aq}{cq}{xq} 
+ (1-c)\, \opo{a}{c}{x}. 
\end{align}
We now prove the uniqueness of $(S, T)$ by contradiction. 
Suppose that two distinct pairs $(S_1,T_1)$ and $(S_2,T_2)$ satisfy \eqref{3trr_1phi1}. 
Then
\[
(S_1-S_2)\, \opo{aq}{cq}{xq} 
=
(T_2-T_1)\, \opo{a}{c}{x}. 
\]
This implies that ${}_1\phi_1 (aq, cq; q, xq) / {}_1\phi_1 (a,c;q,x)$ is a rational function in $a,c,x,$ and $q$. 
Moreover, by \eqref{3trr_1phi1_(1,1,0)}, 
it follows that ${}_1\phi_1 (aq, cq; q, x) / {}_1\phi_1 (a,c;q,x)$ is also a rational function in $a,c,x,$ and $q$. 
However, from the formula (see, e.g., \cite[p.~29, 1.18(i)]{GR})
\begin{align*}
\opo{a}{0}{-q} = \frac{(aq;q^2)_\infty}{(q;q^2)_\infty}, 
\end{align*}
it follows that when $c=0$ and $x=-q$, the ratio ${}_1\phi_1 (aq, cq; q, x) / {}_1\phi_1 (a,c;q,x)$ has an unbounded set of poles as a function of $a$. 
This is a contradiction. 
Hence, the uniqueness of $(S, T)$ satisfying \eqref{3trr_1phi1} is proved. 

Next, we prove the relation \eqref{relation_S,T} between $S$ and $T$. 
To this end, we use the following three-term recurrence relation for ${}_2\phi_1$ obtained from \cite[Theorem~1]{Y1}: 
\begin{align*}
\p{aq^2}{bq^2}{cq^2}{x} 
&= -\frac{(1-cq)\{1-c-bx-a(1-b-bq)x\}}{(1-aq)(1-bq)(c-abxq)x}\, \p{aq}{bq}{cq}{x} \\
&\quad + \frac{(1-c)(1-cq)}{(1-aq)(1-bq)(c-abxq)x}\, \p{a}{b}{c}{x}. 
\end{align*}
Replacing $x$ with $x/b$ and letting $b \to \infty$, we obtain 
\begin{align}\label{3trr_1phi1_(2,2,2)}
\opo{aq^2}{cq^2}{xq^2} 
&= \frac{(1-cq)\{1-c-(1-a-aq)x \}}{(1-aq)(c-axq)xq}\, \opo{aq}{cq}{xq} \\
&\quad - \frac{(1-c)(1-cq)}{(1-aq)(c-axq)xq}\, \opo{a}{c}{x}. \nonumber
\end{align}
We now prove \eqref{relation_S,T}. 
Replacing $(k,m,n,a,c,x)$ in \eqref{3trr_1phi1} with $(k-1,m-1,n-1,aq,cq,xq)$, we obtain 
\begin{align*}
\opo{aq^k}{cq^m}{xq^n}
&= S(k-1,m-1,n-1;aq,cq,xq;q)\, \opo{aq^2}{cq^2}{xq^2} \\
&\quad + T(k-1,m-1,n-1;aq,cq,xq;q)\, \opo{aq}{cq}{xq}. 
\end{align*}
Moreover, by \eqref{3trr_1phi1_(2,2,2)}, we have 
\begin{align*}
\opo{aq^k}{cq^m}{xq^n} 
&= S' \cdot\opo{aq}{cq}{xq} \\
&- \frac{(1-c)(1-cq)}{(1-aq)(c-axq)xq} S(k-1,m-1,n-1;aq,cq,xq;q)\, \opo{a}{c}{x}, \nonumber
\end{align*}
where $S'$ is a certain rational function in $a,c,x,$ and $q$. 
From the uniqueness of $(S, T)$ satisfying \eqref{3trr_1phi1}, 
we obtain \eqref{relation_S,T}. 

Finally, we give the explicit expression \eqref{expression_S} for $S$ by using the results in \cite{Y1}. 
We denote the coefficients of \eqref{3trr_2phi1} by
\[
Q = Q(k,l,m,n; a,b,c,x; q), \qquad
R = R(k,l,m,n; a,b,c,x; q).
\]
Consider the three-term recurrence relation 
\begin{align*}
\p{a q^{k}}{b q^{l'}}{c q^{m}}{x q^{n'}} 
&= Q (k,l',m,n'; a,b,c,x; q)\, \p{a q}{b q}{c q}{x} \\
&\quad + R (k,l',m,n'; a,b,c,x; q)\, \p{a}{b}{c}{x}, 
\end{align*}
where $k,l',m,n'$ are arbitrary integers. 
Replacing $x$ with $x/b$ and letting $b \to \infty$, we obtain 
\begin{align*}
\opo{a q^k}{c q^m}{x q^n}
&= \left(\lim_{b \to \infty}\, Q\!\left(k,l',m,n'; a,b,c,\frac{x}{b}; q\right) \right)\, \opo{a q}{c q}{x q} \\
&\quad 
+ \left(\lim_{b \to \infty}\, R\!\left(k,l',m,n'; a,b,c,\frac{x}{b}; q\right) \right)\, \opo{a}{c}{x}, 
\end{align*}
where $n := l' + n'$. 
By the uniqueness of $(S, T)$ satisfying \eqref{3trr_1phi1}, it follows that 
\begin{align}
S(k,m,n;a,c,x;q) 
&= \lim_{b \to \infty}\, Q\!\left(k,l',m,n'; a,b,c,\frac{x}{b}; q\right), \label{S=limQ}\\
T(k,m,n;a,c,x;q) 
&= \lim_{b \to \infty}\, R\!\left(k,l',m,n'; a,b,c,\frac{x}{b}; q\right). \nonumber
\end{align}
From \cite[Theorem~1]{Y1}, we obtain 
\begin{align}\label{expression_Q}
Q\!\left(k,l',m,n'; a,b,c,\frac{x}{b}; q \right)
&= \frac{(1-a)(1-\frac{1}{b})(c-ax)}{(q-c)(1-c)} 
\frac{(x/b;q)_{\min\{n',0\}} x^{1-\max\{m,0 \}} }{(ax/c;q)_{\max\{k+l'-m+n',0\}}} \\
&\quad \times b^{\max\{m,0\}} P \c{k}{l'}{m}{n'}{a}{b}{c}{\frac{x}{b}}. \nonumber
\end{align}
Here we remark that the condition $k \leq l$ imposed in \cite[Theorem~1]{Y1} can be disregarded. 
The reason is that, in the present proof, we use the explicit expression for the polynomial $P$ in $x$ given in \cite[Lemma~15]{Y1}, 
and this expression is symmetric under the interchange $(a,k) \leftrightarrow (b,l)$.

From \eqref{S=limQ} and \eqref{expression_Q}, it remains only to show that 
\begin{align}\label{limbP=P^(1,1)}
\lim_{b \to \infty} b^{\max\{m,0\}} P \c{k}{l'}{m}{n'}{a}{b}{c}{\frac{x}{b}} 
= P^{(1,1)}(k,m,n;a,c,x;q)
\end{align}
in order to complete the proof of \eqref{expression_S}. 
From \cite[Lemma~15]{Y1}, we obtain 
\begin{align*}
&\lim_{b \to \infty} b^{\max\{m,0\}} P \c{k}{l'}{m}{n'}{a}{b}{c}{\frac{x}{b}} \\
&\quad = x^{\max\{m,0\}} 
\begin{cases}
\displaystyle 
\lim_{b \to \infty} \left(\sum_{j=0}^{\infty} A'_j \left(\frac{x}{b} \right)^{j-m} - \sum_{j=0}^{\infty} B'_j \left(\frac{x}{b} \right)^j\right), 
& k+l'-m+n' \geq 0, \\
\displaystyle 
\lim_{b \to \infty} \left(\sum_{j=0}^{\infty} \tilde{A}'_j \left(\frac{x}{b} \right)^{j-m} - \sum_{j=0}^{\infty} \tilde{B}'_j \left(\frac{x}{b} \right)^j\right), 
& k+l'-m+n' < 0, 
\end{cases} 
\end{align*}
where $A'_j$, $B'_j$, $\tilde{A}'_j$, and $\tilde{B}'_j$ are obtained from
$A_j$, $B_j$, $\tilde{A}_j$, and $\tilde{B}_j$, respectively,
by replacing $(l,n)$ with $(l',n')$. 
The definitions of $A_j$, $B_j$, $\tilde{A}_j$, and $\tilde{B}_j$ are given in \cite[Theorem~1]{Y1}. 
By direct calculation, we obtain
\begin{align*}
\lim_{b\to \infty} b^{m-j} A'_j &= A^{(1,1)}_j, &
\lim_{b\to \infty} b^{-j} B'_j &= B^{(1,1)}_j, \\
\lim_{b\to \infty} b^{m-j} \tilde{A}'_j &= \tilde{A}^{(1,1)}_j, &
\lim_{b\to \infty} b^{-j} \tilde{B}'_j &= \tilde{B}^{(1,1)}_j.
\end{align*}
Hence, we obtain
\begin{align*}
&\lim_{b \to \infty} b^{\max\{m,0\}} P \c{k}{l'}{m}{n'}{a}{b}{c}{\frac{x}{b}} \\
&\quad = x^{\max\{m,0\}} 
\begin{cases}
\displaystyle 
\sum_{j=0}^{\infty} A^{(1,1)}_j x^{j-m} - \sum_{j=0}^{\infty} B^{(1,1)}_j x^j, 
& k-m+n \geq 0, \\
\displaystyle 
\sum_{j=0}^{\infty} \tilde{A}^{(1,1)}_j x^{j-m} - \sum_{j=0}^{\infty} \tilde{B}^{(1,1)}_j x^j, 
& k-m+n < 0.  
\end{cases} \nonumber
\end{align*}
Moreover, by Remark~\ref{remark:P^(1,1)}, it follows that 
\begin{align}\label{limbP}
&\lim_{b \to \infty} b^{\max\{m,0\}} P \c{k}{l'}{m}{n'}{a}{b}{c}{\frac{x}{b}} \\
&\quad = x^{\max\{m,0\}} 
\begin{cases}
\displaystyle 
\sum_{j=-\max\{m,0\}}^{\infty} \left(A^{(1,1)}_{j+m} - B^{(1,1)}_j \right) x^j, 
& k-m+n \geq 0, \\
\displaystyle 
\sum_{j=-\max\{m,0\}}^{\infty} \left(\tilde{A}^{(1,1)}_{j+m} - \tilde{B}^{(1,1)}_j \right) x^j, 
& k-m+n < 0.  
\end{cases} \nonumber
\end{align}
To show that the right-hand side coincides with
$P^{(1,1)}(k,m,n; a,c,x; q)$, we use the following lemma from
\cite[Lemma~4~(1.5)]{Y1}. 
\begin{lemma}\label{Y1_lem4}
Assume that $j \in \mathbb{Z}$ and $r, n_{1}, \dotsc, n_{r}, s \in \mathbb{Z}_{\geq 0}$. 
If \,$\big\lvert b_1^{-1} q^{j + 1 - (n_{1} + \dotsm + n_{r}) - s} \big\rvert < 1$, then 
\begin{align*}
&{}_{r + 2}\phi_{r + 1} \biggl(\genfrac..{0pt}{}{b_1, \,b_2, \,c_{1} q^{n_{1}}, \dotsc, c_{r} q^{n_{r}}}{b_2 q^{j + 1}, \,c_{1}, \dotsc, c_{r}}; q,\, b_1^{-1} q^{j + 1 - (n_{1} + \dotsm + n_{r}) - s}\biggr) \\
&= \frac{(q;q)_{\infty} (b_2 q / b_1;q)_{\infty} (b_2 q;q)_{j} (c_{1} / b_2;q)_{n_{1}} \dotsm (c_{r} / b_2;q)_{n_{r}}}{(q / b_1;q)_{\infty} (b_2 q;q)_{\infty} (q;q)_{j} (c_{1};q)_{n_{1}} \dotsm (c_{r};q)_{n_{r}}} 
b_2^{n_{1} + \dotsm + n_{r} - j + s}  \\
&\quad \times {}_{r + 2}\phi_{r + 1}\biggl(\genfrac..{0pt}{}{q^{-j}, \,b_2, \,b_2 q / c_{1}, \dotsc, b_2 q / c_{r}}{b_2 q / b_1, \,b_2 q^{1 - n_{1}} / c_{1}, \dotsc, b_2 q^{1 - n_{r}} / c_{r}}; q,\, q^{s+1}\biggr). 
\end{align*}
\end{lemma}

Assume that $k-m+n \geq 0$, $j \geq k-m$, and $j \geq n-m$.
Setting $r=1$, $b_1=q^{-m-j}$, $b_2=cq^{-j-1}$, $c_1=aq^{k-m-j}$,
$n_1=j+m-k$, and $s=k-m+n$ in Lemma~\ref{Y1_lem4}, 
and after a straightforward simplification, we obtain
\begin{align*}
&\tpt{q^{-m-j}}{cq^{-j-1}}{a}{c}{aq^{k-m-j}}{q^{j+m-n+1}} \\
&\quad
= \frac{(q;q)_{j+m} (a;q)_{k-m-j} (aq/c;q)_j}
{(c;q)_m (q;q)_j (aq/c;q)_{k-m}}
(cq^{-j-1})^{n}
\tpt{q^{-j}}{cq^{-j-1}}{cq^{m-k}/a}{cq^{m}}{cq^{-j}/a}{q^{k-m+n+1}}.
\end{align*}
This shows that $A^{(1,1)}_{j+m} = B^{(1,1)}_j$
when $k-m+n \geq 0$ and $j \geq \max\{k-m,n-m\}$. 
Hence, for $k-m+n \geq 0$, we obtain 
\begin{align}\label{limbP_1}
&x^{\max\{m,0\}} \sum_{j=-\max\{m,0\}}^{\infty} 
\left(A^{(1,1)}_{j+m} - B^{(1,1)}_j \right) x^j \\
&\quad = \sum_{j=-\max\{m,0\}}^{\max\{k-m,\, n-m\} - 1} 
\left(A^{(1,1)}_{j+m} - B^{(1,1)}_j \right) x^{j+\max\{m,0\}} \nonumber \\
&\quad = \sum_{j=0}^{\max\{k-m,\, n-m\}+\max\{m,0\}-1} 
\left(A^{(1,1)}_{j+\min\{m,0\}} - B^{(1,1)}_{j-\max\{m,0\}} \right) x^{j}. \nonumber
\end{align}
Here we note that 
\[
\max\{k-m,\, n-m\}+\max\{m,0\}-1 = k-m+n-\min\{k,n\}+\max\{m,0\}-1 = d. 
\]
Next, assume that $k-m+n \leq 0$, $j \geq -k$, and $j \geq -n$.
Setting $r=1$, $b_1=q^{-m-j}$, $b_2=cq^{-j-1}$, $c_1=cq^{-k-j}/a$,
$n_1=j+k$, and $s=j+n$ in Lemma~\ref{Y1_lem4}, 
and after a straightforward simplification, we obtain
\begin{align*}
&\tpt{q^{-m-j}}{cq^{-j-1}}{c/a}{c}{cq^{-k-j}/a}{q^{m-k-n+1}} \\
&\quad
= \frac{(q;q)_{j+m} (a;q)_{k} (q/a;q)_j}
{(c;q)_m (q;q)_j (aq/c;q)_{j+k}}
(-1)^j a^{j} c^{n}
q^{-\binom{j+1}{2}-(j+1)n} 
\tpt{q^{-j}}{cq^{-j-1}}{aq^{k}}{cq^{m}}{aq^{-j}}{q^{j+n+1}}.
\end{align*}
This shows that $\tilde{A}^{(1,1)}_{j+m} = \tilde{B}^{(1,1)}_j$
when $k-m+n < 0$ and $j \geq \max\{-k,-n\}$. 
Hence, for $k-m+n < 0$, we obtain 
\begin{align}\label{limbP_2}
&x^{\max\{m,0\}} \sum_{j=-\max\{m,0\}}^{\infty} 
\left(\tilde{A}^{(1,1)}_{j+m} - \tilde{B}^{(1,1)}_j \right) x^j \\
&\quad = \sum_{j=-\max\{m,0\}}^{\max\{-k,\, -n\} - 1} 
\left(\tilde{A}^{(1,1)}_{j+m} - \tilde{B}^{(1,1)}_j \right) x^{j+\max\{m,0\}} \nonumber \\
&\quad = \sum_{j=0}^{\max\{-k,\, -n\}+\max\{m,0\}-1} 
\left(\tilde{A}^{(1,1)}_{j+\min\{m,0\}} - \tilde{B}^{(1,1)}_{j-\max\{m,0\}} \right) x^{j}. \nonumber
\end{align}
Here we note that 
\[
\max\{-k,\, -n\}+\max\{m,0\}-1 = -\min\{k,n\}+\max\{m,0\}-1 = d. 
\]
From \eqref{P^(1,1)} and \eqref{limbP}--\eqref{limbP_2}, we obtain \eqref{limbP=P^(1,1)}. 
This completes the proof of \eqref{expression_S}.

\section{An alternative recurrence relation for ${}_1\phi_1$}
In this section, we prove Theorem~\ref{thm:2}. 

First, we prove by contradiction the uniqueness of the pair $(\tilde{S}, \tilde{T})$ satisfying \eqref{3trr_1phi1_2}. 
Suppose that two distinct pairs $(\tilde{S}_1,\tilde{T}_1)$ and $(\tilde{S}_2,\tilde{T}_2)$ satisfy \eqref{3trr_1phi1_2}. 
Then it follows that ${}_1\phi_1 (aq, cq; q, x) / {}_1\phi_1 (a,c;q,x)$ is a rational function in $a,c,x,$ and $q$. 
However, as shown in Section~$\ref{proof_thm1}$, the ratio has an unbounded set of poles as a function of $a$. 
This is a contradiction. 
Hence, the uniqueness of $(\tilde{S}, \tilde{T})$ satisfying \eqref{3trr_1phi1_2} is proved. 

Next, we prove the relation \eqref{relation_tS,tT} between $\tilde{S}$ and $\tilde{T}$.
To this end, we use the explicit expressions for the coefficients of \eqref{3trr_1phi1} given in Theorem~\ref{thm:1}. 
Applying Theorem~\ref{thm:1} with $(k,m,n)=(1,1,0)$, we obtain
\begin{align}\label{3trr_1phi1_(1,1,0)again}
\opo{aq}{cq}{x}
&= (c-ax)\,\opo{aq}{cq}{xq}
+ (1-c)\,\opo{a}{c}{x}, 
\end{align}
which coincides with \eqref{3trr_1phi1_(1,1,0)} derived earlier. 
Applying Theorem~\ref{thm:1} with $(k,m,n)=(2,2,0)$, we further obtain
\begin{align}\label{3trr_1phi1_(2,2,0)}
\opo{aq^2}{cq^2}{x}
&= \frac{(1-cq)(c-ax)\{(1-c)q+x\}}{(1-aq)x}\,\opo{aq}{cq}{xq} \\
&\quad - \frac{(1-c)(1-cq)(cq-x)}{(1-aq)x}\,\opo{a}{c}{x}. \nonumber
\end{align}
Eliminating the term ${}_1\phi_1(aq,cq;q,xq)$ from these equations, 
we obtain
\begin{align}\label{3trr_1phi1_(2,2,0)_tilde}
\opo{aq^2}{cq^2}{x} 
&= \frac{(1-cq)\{(1-c)q+x \}}{(1-aq)x}\, \opo{aq}{cq}{x} \\
&\quad - \frac{(1-c)(1-cq)q}{(1-aq)x}\, \opo{a}{c}{x}. \nonumber
\end{align}
We now prove \eqref{relation_tS,tT}. 
Replacing $(k,m,a,c)$ in \eqref{3trr_1phi1_2} with $(k-1,m-1,aq,cq)$, we obtain 
\begin{align*}
\opo{aq^k}{cq^m}{xq^n}
&= \tilde{S}(k-1,m-1,n;aq,cq,x;q)\, \opo{aq^2}{cq^2}{x} \\
&\quad + \tilde{T}(k-1,m-1,n;aq,cq,x;q)\, \opo{aq}{cq}{x}. 
\end{align*}
Moreover, by \eqref{3trr_1phi1_(2,2,0)_tilde}, we have 
\begin{align*}
\opo{aq^k}{cq^m}{xq^n} 
&= \tilde{S}' \cdot\opo{aq}{cq}{x} \\
&- \frac{(1-c)(1-cq)q}{(1-aq)x}\, \tilde{S}(k-1,m-1,n;aq,cq,x;q)\, \opo{a}{c}{x}, \nonumber
\end{align*}
where $\tilde{S}'$ is a certain rational function in $a,c,x,$ and $q$. 
From the uniqueness of $(\tilde{S}, \tilde{T})$ satisfying \eqref{3trr_1phi1_2}, 
we obtain \eqref{relation_tS,tT}. 

Finally, we give the explicit expression \eqref{expression_tS} for $\tilde{S}$. 
Using \eqref{3trr_1phi1_2} together with \eqref{3trr_1phi1_(1,1,0)again}, we obtain 
\begin{align*}
\opo{a q^k}{c q^m}{x q^n}
&= \tilde{S} \cdot\opo{a q}{c q}{x}
+ \tilde{T} \cdot\opo{a}{c}{x} \\
&= (c-ax)\, \tilde{S} \; \opo{a q}{c q}{x q} 
+ \{(1-c)\, \tilde{S} + \tilde{T} \}\, \opo{a}{c}{x}. 
\end{align*}
From the uniqueness of $(S, T)$ satisfying \eqref{3trr_1phi1}, 
it follows that 
\begin{align*}
(c-ax)\, \tilde{S}(k,m,n;a,c,x;q) = S(k,m,n;a,c,x;q). 
\end{align*}
Combining this with \eqref{expression_S}, we obtain \eqref{expression_tS}.

\section{Three-term recurrence relation for ${}_0\phi_1$}
In this section, we prove Theorem~\ref{thm:3}. 

First, we prove the uniqueness of the pair $(U, V)$ satisfying \eqref{3trr_0phi1}. 
To prove this, we use the three-term recurrence relation 
\begin{align}\label{3trr_0phi1_(1,2)}
\zpo{cq}{xq} = (c-x)\, \zpo{cq}{xq^2} + (1-c)\, \zpo{c}{x}
\end{align}
obtained by replacing $x$ in \eqref{3trr_1phi1_(1,1,0)again} with $x/a$ and letting $a \to \infty$. 
We now prove the uniqueness of $(U, V)$ by contradiction. 
Suppose that two distinct pairs $(U_1,V_1)$ and $(U_2,V_2)$ satisfy \eqref{3trr_0phi1}. 
Then
\[
(U_1-U_2)\, \zpo{cq}{xq} 
=
(V_2-V_1)\, \zpo{c}{x}. 
\]
This implies that ${}_0\phi_1 (\text{--},c;q,x) / {}_0\phi_1 (\text{--},cq;q,xq)$ is a rational function in $c,x,$ and $q$. 
Moreover, by \eqref{3trr_0phi1_(1,2)}, 
it follows that ${}_0\phi_1 (\text{--},cq;q,xq^2) / {}_0\phi_1 (\text{--},cq;q,xq)$ is also a rational function in $c,x,$ and $q$. 
However, from the formula (see, e.g., \cite[(1.3.16)]{GR})
\begin{align*}
\zpo{-q}{x} 
= {}_{0}\phi_{0}\biggl(\genfrac..{0pt}{}{\text{--}}{\text{--}};q^2, -x\biggr)
= (-x;q^2)_\infty, 
\end{align*}
it follows that when $c=-1$, the ratio ${}_0\phi_1 (\text{--},cq;q,xq^2) / {}_0\phi_1 (\text{--},cq;q,xq)$ has an unbounded set of poles as a function of $x$. 
This is a contradiction. 
Hence, the uniqueness of $(U, V)$ satisfying \eqref{3trr_0phi1} is proved. 

Next, we prove the relation \eqref{relation_U,V} between $U$ and $V$. 
To this end, we use the three-term recurrence relation 
\begin{align}\label{3trr_0phi1_(2,2)}
\zpo{cq^2}{xq^2} = -\frac{(1-c)(1-cq)}{x}\, \zpo{cq}{xq} + \frac{(1-c)(1-cq)}{x}\, \zpo{c}{x}
\end{align}
obtained by replacing $x$ in \eqref{3trr_1phi1_(2,2,0)_tilde} with $x/a$ and letting $a \to \infty$. 
We now prove \eqref{relation_U,V}. 
Replacing $(m,n,c,x)$ in \eqref{3trr_0phi1} with $(m-1,n-1,cq,xq)$, we obtain 
\begin{align*}
\zpo{c q^m}{x q^n}
&= U(m-1,n-1;cq,xq;q) \,\zpo{c q^2}{x q^2} \\
&\quad + V(m-1,n-1;cq,xq;q) \,\zpo{cq}{xq}.
\end{align*}
Moreover, by \eqref{3trr_0phi1_(2,2)}, we have 
\begin{align*}
\zpo{c q^m}{x q^n}
&= U' \cdot\zpo{c q}{x q} \\
&\quad + \frac{(1-c)(1-cq)}{x} U(m-1,n-1;cq,xq;q)\, \zpo{c}{x}, 
\end{align*}
where $U'$ is a certain rational function in $c,x,$ and $q$. 
From the uniqueness of $(U, V)$ satisfying \eqref{3trr_0phi1}, we obtain \eqref{relation_U,V}. 

Finally, we derive the explicit expression \eqref{expression_U} for $U$ as a confluent limit of the explicit expression \eqref{expression_tS} for $\tilde{S}$. 
Consider the three-term recurrence relation 
\begin{align*}
\opo{a q^{k'}}{c q^m}{x q^{n'}}
&= \tilde{S}(k',m,n';a,c,x;q)\, \opo{a q}{c q}{x} \\
&\quad + \tilde{T}(k',m,n';a,c,x;q)\, \opo{a}{c}{x}, 
\end{align*}
where $k',m,n'$ are arbitrary integers. 
Replacing $x$ with $x/a$ and letting $a \to \infty$, we obtain 
\begin{align*}
\zpo{c q^m}{x q^n}
&= \left(\lim_{a \to \infty} \tilde{S}\!\left(k',m,n'; a,c,\frac{x}{a}; q\right) \right)\, \zpo{c q}{x q} \\
&\quad 
+ \left(\lim_{a \to \infty} \tilde{T}\!\left(k',m,n'; a,c,\frac{x}{a}; q\right) \right)\, \zpo{c}{x}, 
\end{align*}
where $n := k' + n'$. 
By the uniqueness of $(U, V)$ satisfying \eqref{3trr_0phi1}, it follows that 
\begin{align}
U(m,n;c,x;q) 
&= \lim_{a\to \infty} \tilde{S}\!\left(k',m,n'; a,c,\frac{x}{a}; q\right), \label{U=limtS}\\
V(m,n;c,x;q) 
&= \lim_{b \to \infty} \tilde{T}\!\left(k',m,n'; a,c,\frac{x}{a}; q\right). \nonumber
\end{align}
From \eqref{expression_tS}, we obtain 
\begin{align}\label{limtS}
\lim_{a\to \infty}\, \tilde{S}\!\left(k',m,n'; a,c,\frac{x}{a}; q\right) 
&= -\frac{1}{(q-c)(1-c)} 
\frac{x^{1-\max\{m, 0\}}}{(x/c;q)_{\max\{k'-m+n', 0\}}} \\
&\quad \times 
\lim_{a\to \infty} a^{\max\{m,0\}} P^{(1,1)}\!\left(k',m,n'; a,c,\frac{x}{a}; q \right). \nonumber
\end{align}
From \eqref{U=limtS} and \eqref{limtS}, it remains only to show that 
\begin{align}\label{limaP^(1,1)=P^(0,1)}
\lim_{a\to \infty} a^{\max\{m,0\}} 
P^{(1,1)}\!\left(k',m,n'; a,c,\frac{x}{a}; q \right) 
= P^{(0,1)}(m,n; c,x; q) 
\end{align}
in order to complete the proof of \eqref{expression_U}. 
By direct calculation, we obtain 
\begin{align*}
\lim_{a\to \infty} a^{m-j} A'^{(1,1)}_j &= A^{(0,1)}_j, &
\lim_{a\to \infty} a^{-j} B'^{(1,1)}_j &= B^{(0,1)}_j, \\
\lim_{a\to \infty} a^{m-j} \tilde{A}'^{(1,1)}_j &= \tilde{A}^{(0,1)}_j, &
\lim_{a\to \infty} a^{-j} \tilde{B}'^{(1,1)}_j &= \tilde{B}^{(0,1)}_j, 
\end{align*}
where $A'^{(1,1)}_j$, $B'^{(1,1)}_j$, $\tilde{A}'^{(1,1)}_j$, and $\tilde{B}'^{(1,1)}_j$ are obtained from
$A^{(1,1)}_j$, $B^{(1,1)}_j$, $\tilde{A}^{(1,1)}_j$, and $\tilde{B}^{(1,1)}_j$, respectively,
by replacing $(k,n)$ with $(k',n')$. 
Hence, by \eqref{P^(1,1)}, when $k'-m+n' \geq 0$, we obtain
\begin{align*}
&\lim_{a\to \infty} a^{\max\{m,0\}} P^{(1,1)}\!\left(k',m,n'; a,c,\frac{x}{a}; q \right) \\
&\quad = 
\lim_{a\to \infty}\sum_{j=\min\{m,0\}}^{d'+\min\{m,0\}} a^{m-j} A'^{(1,1)}_j x^{j-\min\{m,0\}} 
- \lim_{a\to \infty}\sum_{j=-\max\{m,0\}}^{d'-\max\{m,0\}} a^{-j} B'^{(1,1)}_j x^{j+\max\{m,0\}} \\
&\quad = 
\sum_{j=\min\{m,0\}}^{d'+\min\{m,0\}} A^{(0,1)}_j x^{j-\min\{m,0\}} - \sum_{j=-\max\{m,0\}}^{d'-\max\{m,0\}} B^{(0,1)}_j x^{j+\max\{m,0\}} \\ 
&\quad = 
\sum_{j=\min\{m,0\}-m}^{d'+\min\{m,0\}-m} A^{(0,1)}_{j+m} x^{j+m-\min\{m,0\}} - \sum_{j=-\max\{m,0\}}^{d'-\max\{m,0\}} B^{(0,1)}_j x^{j+\max\{m,0\}} \\
&\quad = 
\sum_{j=-\max\{m,0\}}^{d'-\max\{m,0\}} \left(A^{(0,1)}_{j+m} - B^{(0,1)}_j \right) x^{j+\max\{m,0\}}. 
\end{align*}
Here 
\[
d' := \max\{k'-m+n',0\} + \max\{m,0\}-\min\{k',n'\}-1. 
\]
A similar argument applies when $k'-m+n'<0$.
Consequently, we obtain 
\begin{align}\label{limaP^(1,1)}
&\lim_{a\to \infty} a^{\max\{m,0\}} P^{(1,1)}\!\left(k',m,n'; a,c,\frac{x}{a}; q \right) \\
&=\quad 
\begin{cases}
\displaystyle 
\sum_{j=-\max\{m,0\}}^{d'-\max\{m,0\}} \left(A^{(0,1)}_{j+m} - B^{(0,1)}_j \right) x^{j+\max\{m,0\}}, & n-m \geq 0, \\
\displaystyle 
\sum_{j=-\max\{m,0\}}^{d'-\max\{m,0\}} \left(\tilde{A}^{(0,1)}_{j+m} - \tilde{B}^{(0,1)}_j \right) x^{j+\max\{m,0\}}, & n-m<0. 
\end{cases} \nonumber 
\end{align}
To show that the right-hand side coincides with
$P^{(0,1)}(m,n; c,x; q)$, we use Lemma~\ref{Y1_lem4}. 
Assume that $n-m \geq 0$ and $2j \geq n-2m$.
Setting $r=0$, $b_1=q^{-m-j}$, $b_2=cq^{-j-1}$, and $s=n-m$ in Lemma~\ref{Y1_lem4}, 
and after a straightforward simplification, we obtain
\begin{align*}
\p{q^{-m-j}}{cq^{-j-1}}{c}{q^{2j+2m-n+1}} 
= \frac{(q;q)_{j+m}}{(c;q)_m (q;q)_j}(cq^{-j-1})^{n-m-j}
\p{q^{-j}}{cq^{-j-1}}{cq^{m}}{q^{n-m+1}}.
\end{align*}
This shows that $A^{(0,1)}_{j+m} = B^{(0,1)}_j$
when $n-m \geq 0$ and $j \geq \left\lfloor \frac{n+1}{2} \right\rfloor-m$. 
Hence, for $n-m \geq 0$, we obtain
\begin{align}\label{n-m>=0}
&\sum_{j=-\max\{m,0\}}^{d'-\max\{m,0\}} 
\left(A^{(0,1)}_{j+m} - B^{(0,1)}_j \right) x^{j+\max\{m,0\}} \\
&\quad 
= \sum_{j=-\max\{m,0\}}^{\min\{d'-\max\{m,0\},\, \lfloor (n+1)/2 \rfloor-m-1\}} 
\left(A^{(0,1)}_{j+m} - B^{(0,1)}_j \right) x^{j+\max\{m,0\}} \nonumber \\
&\quad 
= \sum_{j=0}^{\min\{d',\, \lfloor (n+1)/2 \rfloor-m+\max\{m,0\}-1\}} 
\left(A^{(0,1)}_{j+m-\max\{m,0\}} - B^{(0,1)}_{j-\max\{m,0\}} \right) x^{j}. \nonumber
\end{align}
Since 
\begin{align*}
d' 
= (k' + n' - \min\{k',n'\}) -m + \max\{m,0\} - 1 
\geq 
\left\lfloor \frac{n+1}{2} \right\rfloor -m+\max\{m,0\}-1 
\end{align*}
for $n-m \geq 0$, we note that 
\begin{align*}
\min\left\{d', \left\lfloor \frac{n+1}{2} \right\rfloor-m+\max\{m,0\}-1 \right\} 
= \left\lfloor \frac{n+1}{2} \right\rfloor-m+\max\{m,0\}-1 
= e. 
\end{align*}
Next, assume that $n-m \leq 0$ and $2j \geq -n$.
Setting $r=0$, $b_1=q^{-m-j}$, $b_2=cq^{-j-1}$, and $s=2j+n$ in Lemma~\ref{Y1_lem4}, 
and after a straightforward simplification, we obtain
\begin{align*}
\p{q^{-m-j}}{cq^{-j-1}}{c}{q^{m-n+1}} 
= \frac{(q;q)_{j+m}}{(c;q)_m (q;q)_j}(cq^{-j-1})^{j+n}
\p{q^{-j}}{cq^{-j-1}}{cq^{m}}{q^{2j+n+1}}.
\end{align*}
This shows that $\tilde{A}^{(0,1)}_{j+m} = \tilde{B}^{(0,1)}_j$
when $n-m < 0$ and $j \geq \left\lfloor \frac{1-n}{2} \right\rfloor$. 
Hence, for $n-m < 0$, we obtain
\begin{align}\label{n-m<0}
&\sum_{j=-\max\{m,0\}}^{d'-\max\{m,0\}} 
\left(\tilde{A}^{(0,1)}_{j+m} - \tilde{B}^{(0,1)}_j \right) x^{j+\max\{m,0\}} \\
&\quad 
= \sum_{j=-\max\{m,0\}}^{\min\{d'-\max\{m,0\},\, \lfloor (1-n)/2 \rfloor-1\}} 
\left(\tilde{A}^{(0,1)}_{j+m} - \tilde{B}^{(0,1)}_j \right) x^{j+\max\{m,0\}} \nonumber \\
&\quad 
= \sum_{j=0}^{\min\{d',\, \lfloor (1-n)/2 \rfloor+\max\{m,0\}-1\}} 
\left(\tilde{A}^{(0,1)}_{j+m-\max\{m,0\}} - \tilde{B}^{(0,1)}_{j-\max\{m,0\}} \right) x^{j}. \nonumber
\end{align}
Since 
\begin{align*}
d' 
= - \min\{k',n'\} + \max\{m,0\} - 1 \geq 
\left\lfloor \frac{1-n}{2} \right\rfloor +\max\{m,0\}-1 
\end{align*}
for $n-m < 0$, we note that 
\begin{align*}
\min\left\{d', \left\lfloor \frac{1-n}{2} \right\rfloor+\max\{m,0\}-1 \right\} 
= \left\lfloor \frac{1-n}{2} \right\rfloor+\max\{m,0\}-1 
= e. 
\end{align*}
From \eqref{P^(0,1)} and \eqref{limaP^(1,1)}--\eqref{n-m<0}, we obtain \eqref{limaP^(1,1)=P^(0,1)}. 
This completes the proof of \eqref{expression_U}.

\section{Three-term recurrence relation for $J^{(3)}_\nu$}
In this section, we derive Theorem~\ref{thm:4}. 

Substituting $a=0$, $c=q^\nu$, and replacing $x$ with $x^2q/4$ in  \eqref{3trr_1phi1_2}, we obtain 
\begin{align}\label{3trr_J^(3)_pre}
\opo{0}{q^{\nu+m}}{\frac{x^2q^{n+1}}{4}}
= \tilde{s}\cdot \opo{0}{q^{\nu+1}}{\frac{x^2q}{4}} 
+ \tilde{t}\cdot \opo{0}{q^\nu}{\frac{x^2q}{4}}, 
\end{align}
where 
\begin{align*}
\tilde{s} := \tilde{S}\!\left(k,m,n; 0,q^\nu,\frac{x^2q}{4}; q\right), \qquad
\tilde{t} := \tilde{T}\!\left(k,m,n; 0,q^\nu,\frac{x^2q}{4}; q\right). 
\end{align*}

\begin{remark}
The singularity at $a=0$ appearing in the coefficients of \eqref{3trr_1phi1_2} 
(more precisely, in $B^{(1,1)}_j$ and $\tilde{A}^{(1,1)}_j$) is removable,
since the singularity due to $1/a$ cancels in the terminating basic hypergeometric series. 
Hence, the coefficients admit well-defined limits as $a\to0$, and we extend the coefficients to $a=0$ by these limits. 
\end{remark}

As will be seen later, $\tilde{s}$ does not depend on $k$. Hence we write
\[
\tilde{s}=\tilde{s}_{m,n,\nu}(x;q).
\]
Then, by the relation \eqref{relation_tS,tT} between $\tilde{S}$ and $\tilde{T}$, we have
\begin{align*}
\tilde{t}
&= -(q^{\nu};q)_2 \left(\frac{x}{2}\right)^{-2}
\tilde{S}\!\left(k-1,m-1,n; 0,q^{\nu+1},\frac{x^2q}{4}; q\right) \\
&= -(q^{\nu};q)_2 \left(\frac{x}{2}\right)^{-2}\tilde{s}_{m-1,n,\nu+1}(x;q).
\end{align*}
Hence, multiplying both sides of \eqref{3trr_J^(3)_pre} by 
\[
\frac{(q^{\nu+m};q)_\infty}{(q;q)_\infty}
\left(\frac{x}{2}\right)^{\nu+m-1},
\]
we obtain the following relation: 
\begin{align*}
q^{-\frac{(\nu+m-1)n}{2}} J^{(3)}_{\nu+m-1}(xq^{\frac{n}{2}};q) 
&= \frac{(q^{\nu+m};q)_\infty}{(q^{\nu+1};q)_\infty} \left(\frac{x}{2}\right)^{m-1}\,
\tilde{s}_{m,n,\nu}(x;q)\, J^{(3)}_\nu (x;q) \\
&\quad - \frac{(q^{\nu+m};q)_\infty}{(q^{\nu+2};q)_\infty} \left(\frac{x}{2}\right)^{m-2}\,
\tilde{s}_{m-1,n,\nu+1}(x;q)\, J^{(3)}_{\nu-1}(x;q). 
\end{align*}
Replacing $m$ with $m+1$, we obtain 
\begin{align*}
q^{-\frac{(\nu+m)n}{2}} J^{(3)}_{\nu+m}(xq^{\frac{n}{2}};q) 
&= \frac{(q^{\nu+m+1};q)_\infty}{(q^{\nu+1};q)_\infty} \left(\frac{x}{2}\right)^{m}\,
\tilde{s}_{m+1,n,\nu}(x;q)\, J^{(3)}_\nu (x;q) \\
&\quad - \frac{(q^{\nu+m+1};q)_\infty}{(q^{\nu+2};q)_\infty} \left(\frac{x}{2}\right)^{m-1}\,
\tilde{s}_{m,n,\nu+1}(x;q)\, J^{(3)}_{\nu-1}(x;q). 
\end{align*}
Thus, to complete the proof of Theorem~\ref{thm:4}, 
it remains to show that 
\begin{align}\label{ts=R^(3)}
\frac{(q^{\nu+m+1};q)_\infty}{(q^{\nu+1};q)_\infty} \left(\frac{x}{2}\right)^{m}\,
\tilde{s}_{m+1,n,\nu}(x;q) 
= R^{(3)}_{m,n,\nu}(x; q). 
\end{align}
From \eqref{expression_tS}, we obtain 
\begin{align*}
\tilde{s}_{m,n,\nu}(x;q) 
&= \tilde{S}\!\left(k,m,n; 0,q^\nu,\frac{x^2q}{4}; q\right) \\
&= \frac{q^{-\max\{m,0\}}}{(q^{\nu-1};q)_2} \left(\frac{x}{2}\right)^{2-2\max\{m,0\}} 
P^{(1,1)}\!\left(k,m,n; 0,q^\nu,\frac{x^2q}{4}; q\right). 
\end{align*}
By direct calculation, when $c= q^\nu$, we obtain 
\begin{align*}
\lim_{a\to 0} A^{(1,1)}_j &= \lim_{a\to 0} \tilde{A}^{(1,1)}_j = A^{(3)}_{m,n,\nu,j}, \\
\lim_{a\to 0} B^{(1,1)}_j &= \lim_{a\to 0} \tilde{B}^{(1,1)}_j = B^{(3)}_{m,n,\nu,j}.  
\end{align*}
Hence, we obtain 
\begin{align*}
\tilde{s}_{m,n,\nu}(x;q) 
&= \frac{q^{-\max\{m,0\}}}{(q^{\nu-1};q)_2} \left(\frac{x}{2}\right)^{2-2\max\{m,0\}} 
\sum_{j=0}^{d} \left(A^{(3)}_{m,n,\nu,j+\min\{m,0\}} - B^{(3)}_{m,n,\nu,j-\max\{m,0\}}\right) 
\left(\frac{x^2q}{4}\right)^j \\
&= \frac{q^{-\max\{m,0\}}}{(q^{\nu-1};q)_2} \left(\frac{x}{2}\right)^{2-2\max\{m,0\}} 
\sum_{j=-\max\{m,0\}}^{d-\max\{m,0\}} \left(A^{(3)}_{m,n,\nu,j+m} - B^{(3)}_{m,n,\nu,j}\right) 
\left(\frac{x^2q}{4}\right)^{j+\max\{m,0\}}, 
\end{align*}
where $d = \max\{k-m+n, 0\} + \max\{m,0\} - \min\{k,n\} - 1$. 
To show that the sum on the right-hand side coincides with $P^{(3)}_{m,n,\nu}(x;q)$,
we use Lemma~\ref{Y1_lem4}. 
Assume that $j \geq n-m$ and $j \geq -n$.
Setting $r=0$, $b_1=q^{-m-j}$, $b_2=q^{\nu-j-1}$, and $s=j+n$ in Lemma~\ref{Y1_lem4}, 
and after a straightforward simplification, we obtain
\begin{align*}
\p{q^{-m-j}}{q^{\nu-j-1}}{q^\nu}{q^{j+m-n+1}} 
= \frac{(q;q)_{j+m}}{(q^\nu;q)_m (q;q)_j}q^{(\nu-j-1)n}\,
\p{q^{-j}}{q^{\nu-j-1}}{q^{\nu+m}}{q^{j+n+1}}.
\end{align*}
This shows that $A^{(3)}_{m,n,\nu,j+m} = B^{(3)}_{m,n,\nu,j}$ for $j \geq \max\{n-m,-n \}$. 
Since 
\begin{align*}
d-\max\{m,0\} 
&= \max\{k-m+n,0\} + \max\{-k,-n\} - 1 \\
&= \max\{(k-m+n)-k,(k-m+n)-n,-k,-n\} - 1 \\
&= \max\{n-m,k-m,-k,-n\} - 1 \\
&\geq \max\{n-m,-n\}-1, 
\end{align*}
we obtain 
\begin{align*}
&\sum_{j=-\max\{m,0\}}^{d-\max\{m,0\}} \left(A^{(3)}_{m,n,\nu,j+m} - B^{(3)}_{m,n,\nu,j}\right) 
\left(\frac{x^2q}{4}\right)^{j+\max\{m,0\}} \\
&\quad = 
\sum_{j=-\max\{m,0\}}^{\min\{d-\max\{m,0\},\, \max\{n-m,-n\}-1\}} \left(A^{(3)}_{m,n,\nu,j+m} - B^{(3)}_{m,n,\nu,j}\right) 
\left(\frac{x^2q}{4}\right)^{j+\max\{m,0\}} \\
&\quad = 
\sum_{j=-\max\{m,0\}}^{\max\{n-m,-n\}-1} \left(A^{(3)}_{m,n,\nu,j+m} - B^{(3)}_{m,n,\nu,j}\right) 
\left(\frac{x^2q}{4}\right)^{j+\max\{m,0\}} \\
&\quad = 
\sum_{j=0}^{\max\{n-m,-n\}+\max\{m,0\}-1} \left(A^{(3)}_{m,n,\nu,j+m-\max\{m,0\}} - B^{(3)}_{m,n,\nu,j-\max\{m,0\}}\right) 
\left(\frac{x^2q}{4}\right)^{j} \\
&\quad = P^{(3)}_{m,n,\nu}(x;q). 
\end{align*}
Therefore, we obtain 
\begin{align*}
\tilde{s}_{m,n,\nu}(x;q) 
= \frac{q^{-\max\{m,0\}}}{(q^{\nu-1};q)_2} \left(\frac{x}{2}\right)^{2-2\max\{m,0\}} 
P^{(3)}_{m,n,\nu}(x;q), 
\end{align*}
which proves \eqref{ts=R^(3)} and completes the proof of Theorem~\ref{thm:4}.

\section{Three-term recurrence relation for $J^{(2)}_\nu$}
In this section, we derive Theorem~\ref{thm:5}. 

Substituting $c=q^\nu$, replacing $x$ with $-x^2q^\nu/4$, and shifting $n$ by $m$ in \eqref{3trr_0phi1}, we obtain
\begin{align}\label{3trr_J^(2)_pre}
\zpo{q^{\nu+m}}{-\frac{x^2q^{\nu+m+n}}{4}}
&= u_{m,n,\nu}(x;q)\, \zpo{q^{\nu+1}}{-\frac{x^2q^{\nu+1}}{4}} \\
&\quad + v_{m,n,\nu}(x;q)\, \zpo{q^\nu}{-\frac{x^2q^\nu}{4}}, \nonumber
\end{align}
where 
\begin{align*}
u_{m,n,\nu}(x;q) &:= U\!\left(m,m+n;q^\nu,-\frac{x^2q^{\nu}}{4};q \right), \\
v_{m,n,\nu}(x;q) &:= V\!\left(m,m+n;q^\nu,-\frac{x^2q^{\nu}}{4};q \right). 
\end{align*}
By the relation \eqref{relation_U,V} between $U$ and $V$, we have
\begin{align*}
v_{m,n,\nu}(x;q)
&= - (q^\nu; q)_2\, q^{-\nu} \left(\frac{x}{2}\right)^{-2} 
U\!\left(m-1,m+n-1;q^{\nu+1},-\frac{x^2q^{\nu+1}}{4};q \right) \\
&= -(q^{\nu};q)_2\, q^{-\nu} \left(\frac{x}{2}\right)^{-2} u_{m-1,n,\nu+1}(x;q).
\end{align*}
Hence, multiplying both sides of \eqref{3trr_J^(2)_pre} by 
\[
\frac{(q^{\nu+m};q)_\infty}{(q;q)_\infty} q^{(m-1)\nu + \binom{m-1}{2}} 
\left(\frac{x}{2}\right)^{\nu+m-1},
\]
we obtain 
\begin{align*}
&q^{-\frac{(\nu+m-1)n}{2}+(m-1)\nu+\binom{m-1}{2}} J^{(2)}_{\nu+m-1}(xq^{\frac{n}{2}};q) \\
&\quad = \frac{(q^{\nu+m};q)_\infty}{(q^{\nu+1};q)_\infty} q^{(m-1)\nu+\binom{m-1}{2}} 
\left(\frac{x}{2}\right)^{m-1}\, u_{m,n,\nu}(x;q)\, J^{(2)}_\nu (x;q) \\
&\qquad - \frac{(q^{\nu+m};q)_\infty}{(q^{\nu+2};q)_\infty} q^{(m-2)(\nu+1)+\binom{m-2}{2}} 
\left(\frac{x}{2}\right)^{m-2}\, u_{m-1,n,\nu+1}(x;q)\, J^{(2)}_{\nu-1}(x;q). 
\end{align*}
Replacing $m$ with $m+1$, we obtain 
\begin{align*}
&q^{-\frac{(\nu+m)n}{2}+m\nu+\binom{m}{2}} J^{(2)}_{\nu+m}(xq^{\frac{n}{2}};q) \\
&\quad = \frac{(q^{\nu+m+1};q)_\infty}{(q^{\nu+1};q)_\infty} q^{m\nu+\binom{m}{2}} 
\left(\frac{x}{2}\right)^{m}\, u_{m+1,n,\nu}(x;q)\, J^{(2)}_\nu (x;q) \\
&\qquad - \frac{(q^{\nu+m+1};q)_\infty}{(q^{\nu+2};q)_\infty} q^{(m-1)(\nu+1)+\binom{m-1}{2}} 
\left(\frac{x}{2}\right)^{m-1}\, u_{m,n,\nu+1}(x;q)\, J^{(2)}_{\nu-1}(x;q). 
\end{align*}
Thus, to complete the proof of Theorem~\ref{thm:5}, 
it remains to show that 
\begin{align}\label{u=R^(2)}
\frac{(q^{\nu+m+1};q)_\infty}{(q^{\nu+1};q)_\infty} q^{m\nu+\binom{m}{2}} 
\left(\frac{x}{2}\right)^{m}\, u_{m+1,n,\nu}(x;q)
= R^{(2)}_{m,n,\nu}(x; q). 
\end{align}
From \eqref{expression_U}, we obtain 
\begin{align*}
u_{m,n,\nu}(x;q) 
&= U\!\left(m,m+n;q^\nu,-\frac{x^2q^{\nu}}{4};q \right) \\
&= \frac{1}{(q^{\nu-1};q)_2 (-x^2/4;q)_{\max\{n,0\}}} \\
&\quad \times (-1)^{\max\{m,0\}} q^{(1-\max\{m,0\})\nu-1}
\left(\frac{x}{2}\right)^{2-2\max\{m,0\}} 
P^{(0,1)}\!\left(m,m+n; q^\nu,-\frac{x^2q^{\nu}}{4}; q\right). 
\end{align*}
Moreover, by the definitions of $P^{(0,1)}$ and $P^{(2)}$, we obtain 
\begin{align*}
u_{m,n,\nu}(x;q) 
&= \frac{1}{(q^{\nu-1};q)_2 (-x^2/4;q)_{\max\{n,0\}}} \\
&\quad \times (-1)^{\max\{m,0\}} q^{(1-\max\{m,0\})\nu-1}
\left(\frac{x}{2}\right)^{2-2\max\{m,0\}} 
P^{(2)}_{m,n,\nu}(x;q), 
\end{align*}
which proves \eqref{u=R^(2)} and completes the proof of Theorem~\ref{thm:5}.

\clearpage

\section*{Acknowledgment}
This work was supported by JSPS KAKENHI Grant Number JP25KJ0266. 

\section*{Data Availability Statement}
This paper has no associated data. 

\section*{Competing interests}
The author has no competing interests to declare that are relevant to the content of this paper. 

\bibliographystyle{plain}
\bibliography{reference}

@book {GR,
    AUTHOR = {Gasper, George and Rahman, Mizan},
     TITLE = {Basic hypergeometric series},
    SERIES = {Encyclopedia of Mathematics and its Applications},
    VOLUME = {35},
      NOTE = {With a foreword by Richard Askey},
 PUBLISHER = {Cambridge University Press, Cambridge},
      YEAR = {1990},
     PAGES = {xx+287},
      ISBN = {0-521-35049-2},
   MRCLASS = {33Dxx (05A30 33-01 33-02)},
  MRNUMBER = {1052153},
MRREVIEWER = {W.\ A.\ Al-Salam},
}

@article {Ismail,
    AUTHOR = {Ismail, Mourad E. H.},
     TITLE = {The zeros of basic {B}essel functions, the functions
              {$J\sb{\nu +ax}(x)$}, and associated orthogonal polynomials},
   JOURNAL = {J. Math. Anal. Appl.},
  FJOURNAL = {Journal of Mathematical Analysis and Applications},
    VOLUME = {86},
      YEAR = {1982},
    NUMBER = {1},
     PAGES = {1--19},
      ISSN = {0022-247X},
   MRCLASS = {33A65 (33A40)},
  MRNUMBER = {649849},
MRREVIEWER = {M.\ E.\ Muldoon},
       DOI = {10.1016/0022-247X(82)90248-7},
       URL = {https://doi.org/10.1016/0022-247X(82)90248-7},
}

@book {Ismail_book,
    AUTHOR = {Ismail, Mourad E. H.},
     TITLE = {Classical and quantum orthogonal polynomials in one variable},
    SERIES = {Encyclopedia of Mathematics and its Applications},
    VOLUME = {98},
      NOTE = {With two chapters by Walter Van Assche,
              With a foreword by Richard A. Askey},
 PUBLISHER = {Cambridge University Press, Cambridge},
      YEAR = {2005},
     PAGES = {xviii+706},
      ISBN = {978-0-521-78201-2; 0-521-78201-5},
   MRCLASS = {33-02 (05A30 05E35 33C45 33D50 42C05)},
  MRNUMBER = {2191786},
MRREVIEWER = {Bruce\ C.\ Berndt},
       DOI = {10.1017/CBO9781107325982},
       URL = {https://doi.org/10.1017/CBO9781107325982},
}

@article {KS,
    AUTHOR = {Koelink, H. T. and Swarttouw, R. F.},
     TITLE = {On the zeros of the {H}ahn-{E}xton {$q$}-{B}essel function and
              associated {$q$}-{L}ommel polynomials},
   JOURNAL = {J. Math. Anal. Appl.},
  FJOURNAL = {Journal of Mathematical Analysis and Applications},
    VOLUME = {186},
      YEAR = {1994},
    NUMBER = {3},
     PAGES = {690--710},
      ISSN = {0022-247X,1096-0813},
   MRCLASS = {33D45 (33C45)},
  MRNUMBER = {1293849},
MRREVIEWER = {Wenchang\ Chu},
       DOI = {10.1006/jmaa.1994.1327},
       URL = {https://doi.org/10.1006/jmaa.1994.1327},
}

@article {Y0,
    AUTHOR = {Yamaguchi, Yuka},
     TITLE = {Three-term relations for basic hypergeometric series},
   JOURNAL = {J. Math. Anal. Appl.},
  FJOURNAL = {Journal of Mathematical Analysis and Applications},
    VOLUME = {464},
      YEAR = {2018},
    NUMBER = {1},
     PAGES = {662--678},
      ISSN = {0022-247X,1096-0813},
   MRCLASS = {33C05 (33D15)},
  MRNUMBER = {3794109},
MRREVIEWER = {Renu\ Jain},
       DOI = {10.1016/j.jmaa.2018.04.021},
       URL = {https://doi.org/10.1016/j.jmaa.2018.04.021},
}

@article {Y1,
    AUTHOR = {Yamaguchi, Yuka},
     TITLE = {Transformation formulas and three-term relations for basic hypergeometric series},
   JOURNAL = {Funkcial. Ekvac.},
  FJOURNAL = {Funkcialaj Ekvacioj. Serio Internacia},
    VOLUME = {65},
      YEAR = {2022},
    NUMBER = {1},
     PAGES = {35--61},
      ISSN = {0532-8721},
   MRCLASS = {33D15},
  MRNUMBER = {4451007},
       DOI = {10.1619/fesi.65.35},
       URL = {https://doi.org/10.1619/fesi.65.35},
}

\vspace{1cm}
\begin{flushleft}
Faculty of Education \\ 
University of Miyazaki \\ 
1-1 Gakuen Kibanadai-nishi \\ 
Miyazaki 889-2192 Japan \\ 
{\it Email address}: y-yamaguchi@miyazaki-u.ac.jp 
\end{flushleft}

\end{document}